\documentclass{lms}

\usepackage{amsmath, amssymb, amsfonts}
\usepackage{rotating}
\usepackage{url}
\usepackage{microtype}
\usepackage[usenames,dvipsnames]{color}
\usepackage{sistyle}\SIthousandsep{,}
\usepackage[colorlinks]{hyperref}

\usepackage{float}
\newcommand{\fig}[2]{\centering\vspace{.5ex}
\includegraphics[width=#1\textwidth]{#2}}

\renewcommand{\phi}{\varphi}

\newcommand{\Les}{L_E(s)}

\DeclareMathOperator{\tors}{tors}

\DeclareMathOperator{\sinc}{sinc}
\DeclareMathOperator{\Li}{Li}

\def\Z{{\mathbb Z}}

\def\C{{\mathbb C}}
\def\Q{{\mathbb Q}}
\def\F{{\mathbb F}}
\def\Fs{{\mathcal F}}

\usepackage[OT2,T1]{fontenc}
\DeclareSymbolFont{cyrletters}{OT2}{wncyr}{m}{n}
\DeclareMathSymbol{\Sha}{\mathalpha}{cyrletters}{"58}

\def\rk{{\operatorname{rk}}}

\def\Sel{{S}}

\title[Databases of Elliptic Curves Ordered by Height]{Databases of Elliptic Curves Ordered by Height and Distributions of Selmer Groups and Ranks}
\author[Balakrishnan, Ho, Kaplan, Spicer, Stein, Weigandt]{Jennifer S. Balakrishnan, Wei Ho, Nathan Kaplan, Simon Spicer, William Stein and James Weigandt}

\newtheorem{thm}{Theorem}[section] 
\newtheorem{lemma}[thm]{Lemma}     
\newtheorem{cor}[thm]{Corollary}

\newnumbered{conj}{Conjecture}  
\newnumbered{defn}{Definition}
\newnumbered{hypothesis}{Hypothesis}
\newnumbered{remark}{Remark}
\newnumbered{note}{Note}
\newnumbered{observation}{Observation}
\newnumbered{problem}{Problem}
\newnumbered{question}{Question}
\newnumbered{algorithm}{Algorithm}
\newnumbered{example}{Example}
\newunnumbered{notation}{Notation} 

\classno{11G05 (primary), 11-04 (secondary)}

\extraline{We thank the AMS Mathematics Research Communities program, where we began working on this project, for their continued financial support. The first author was partially supported by NSF grant DMS-1103831.  The second author was supported by NSF grant DMS-1406066.}

\begin{document}
\maketitle

\begin{abstract}
Most systematic tables of data associated to ranks of elliptic curves order the curves by conductor.
Recent developments, led by work of Bhargava--Shankar studying the average sizes of $n$-Selmer groups, have given new upper bounds on the average algebraic rank in families of elliptic curves over $\Q$ ordered by height.  We describe databases of elliptic curves over $\Q$ ordered by height in which we compute ranks and $2$-Selmer group sizes, the distributions of which may also be compared to these theoretical results. A striking new phenomenon observed in these databases is that the average rank eventually decreases as height increases.
\end{abstract}


\section{Introduction and Statement of Main Results} \label{sec:intro}

Over the past several decades, tables of elliptic curves defined over $\Q$ have been very useful in number-theoretic research.  A natural ordering on elliptic curves is given by their conductor. Some of the earliest tables were those in Antwerp IV \cite{Antwerp}, which include all elliptic curves of conductor at most $200$.  In \cite{Cremona-book}, Cremona describes algorithms to list all elliptic curves of given conductor and collect arithmetic data for these curves; these algorithms have now produced an exhaustive list of curves of conductor at most $\num{380000}$ in an ongoing project \cite{Cremona}.
A large currently available database of elliptic curves is due to Stein--Watkins \cite{BMSW,SteinWatkins} and includes $\num{136832795}$ curves over $\Q$ of conductor up to $10^8$ and a table of $\num{11378911}$ elliptic curves over $\Q$ of prime conductor up to $10^{10}$, extending earlier tables of this type by Brumer--McGuinness \cite{BrumerMcGuinness}.

It is computationally difficult to produce exhaustive lists of curves up to a given conductor. It is far easier to produce large tables by writing down elliptic curves in Weierstrass form with relatively small defining coefficients; such tables may not, however, include all curves up to a given conductor, as it occasionally happens that curves with small coefficients can have large conductor, or curves with large coefficients can have cancellation leading to a smaller than expected conductor.  For example, the Stein--Watkins table contains approximately $78.5\%$ of the elliptic curves of conductor up to $\num{120000}$ \cite{BMSW}. In very recent work \cite{BeRe}, Bennett and Rechnitzer pursue a different strategy for producing extensive lists of curves with good reduction outside a given prime $p$ by using the reduction theory of binary cubic forms and solving certain Thue-Mahler equations. They find $\num{435893911}$ isomorphism classes of curves of prime conductor up to $10^{12}$ and explain that it is unlikely that any have been missed.

In this paper, we instead describe databases of curves in families ordered by {\em height}---a measure of the size of the coefficients of the Weierstrass equation defining the curve---since it is possible to list all curves in a specified height range. We first consider the family $\Fs_0$ of all elliptic curves over $\Q$ in short Weierstrass form:
\[
\Fs_0 = \{E : y^2 = x^3 + a_4 x + a_6\ |\ a_4,a_6 \in \Z,\ \Delta_E \neq 0\},
\]
where $\Delta_E = -16(4 a_4^3 + 27 a_6^2)$ denotes the discriminant of the curve $E$.
There are two main height functions that we will consider for this family throughout this paper.  The \emph{(naive) height} is defined by $H(E) := \max\{4|a_4|^3, 27 a_6^2\}$.  The \emph{uncalibrated height} is defined by $\widetilde{H}(E) := \max\{|a_4|^3, a_6^2\}$.  We believe that the naive height is the more natural of these two, but both are used in practice.  Both theoretical and computational results using either height have a very similar form.  Throughout this paper we often just write ``height'' in place of ``naive height.''

The main result of this project is the creation of an exhaustive database of isomorphism classes of elliptic curves with naive height up to $2.7\cdot 10^{10}$, a total of $\num{238764310}$ curves \cite{ecdb}. For each elliptic curve in this database, we have recorded its minimal model, torsion subgroup, conductor, Tamagawa product, rank, and $2$-Selmer group rank. The recorded rank for some of the elliptic curves is conditional on the Generalized Riemann Hypothesis (GRH), the Birch and Swinnerton-Dyer (BSD) Conjecture, and/or the Parity Conjecture (see \S \ref{Computing}); however, the recorded rank is unconditional for at least $\num{192850627}$ (approximately $80.77\%$) of the curves. Note that the largest conductor occurring in our database is $\num{863347196528}$.

\begin{figure}[ht]\label{average_rank_figure}
\caption{Average rank of elliptic curves up to a given height}
\fig{.95}{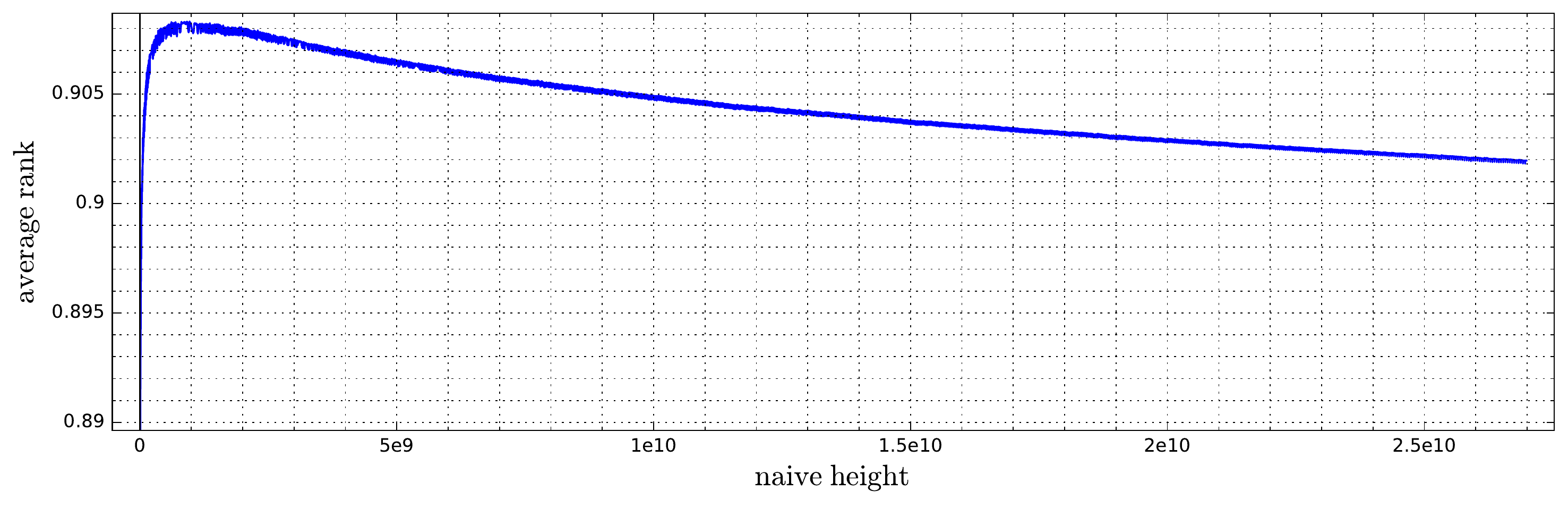}
\end{figure}

In the databases of curves that currently exist, e.g., those compiled in \cite{BMSW,Cremona,SteinWatkins}, the average rank of elliptic curves appeared---from a distance---to be monotonically increasing as the conductor increases, although that would contradict widely believed conjectures (see \S \ref{sec:minimalist}).  Ours is the first database in which we can see a ``turnaround'' point for the average rank of elliptic curves: the average rank of all curves of height up to $X$ appears to be an increasing function of $X$ for $X$ up to approximately $6 \cdot 10^8$ and then looks from a distance like a monotonically decreasing function (though, of course, up close it is wildly gyrating). See Figure \ref{average_rank_figure} for a plot of the average rank of elliptic curves up to height $X$, using our database.  It would be interesting to have theoretical results confirm this observed turnaround point.  Of course, we cannot prove that average rank decreases monotonically (on a large scale) after this turnaround point, but this seems to be a reasonable assumption based on standard heuristics and our data (see \S \ref{sec:minimalist} and \S \ref{sec:samples1}).

From our database, we find that the $\num{238764310}$ curves of naive height up to $2.7 \cdot 10^{10}$ have average rank approximately $0.901976$. The proportion of these curves with each rank is as follows:
\[
{\small
\begin{tabular}{ccccccc}
 Rank $0$ & Rank $1$ & Rank $2$ & Rank $3$ & Rank $4$ & Rank $5$ & Rank $6$ \\
\hline
 $0.32685$ & $0.47381$ & $0.17151$ & $0.02615$ & $0.00159$ & $0.00003$ & $0.000005$
\end{tabular}}\ .
\]
In the appendix, we include additional data about rank distribution for different height ranges. As a point of comparison, the average rank of the $\num{126427408}$ curves with uncalibrated height up to $10^9$ is $0.89473$.

There has been much interest in finding curves of minimal conductor and given rank.  Only the first few of these minimal conductors are known: $11$, $37$, $389$, $5077$, and $234446$.  It is easy to extract analogous results for curves of minimal height from our database.  We include a table of minimal height curves with given torsion subgroup and rank as Table \ref{Records}.

Another problem that arises when computing exhaustive tables of curves is that it is unknown how many curves they will contain if ordered by conductor or discriminant. For example, Watkins \cite[Heuristic 4.1]{Wat} suggests that the number of curves over $\Q$ with conductor bounded by $X$ is asymptotically a constant times $X^{5/6}$, and Brumer--McGuinness \cite{BrumerMcGuinness} conjecture a similar asymptotic for the number of curves over $\Q$ with absolute discriminant less than $X$. In contrast, for curves ordered by height, the number of curves of height at most $X$ is known to be asymptotic to a constant times $X^{5/6}$.  For curves ordered by uncalibrated height, the constant is known to be $4/\zeta(10)$ (see, e.g., \cite{HarSno}). In our database, we observe fast convergence to this asymptotic: the number of curves of uncalibrated height at most $10^9$ is $\num{126427408}$, which is $1.00049\cdot (4/\zeta(10)) (10^9)^{5/6}$.

\subsection{The Minimalist Conjecture} \label{sec:minimalist}

The data sets previously available for elliptic curves defined over $\Q$ are in tension with certain widely believed conjectures (see \cite{BMSW}). We now describe some of these conjectures about the distribution of ranks.

The Minimalist Conjecture, inspired partly by work of Katz--Sarnak relating elliptic curves over function fields to certain random matrix statistics \cite{KatzSarnak}, and by a similar conjecture of Goldfeld for ranks of elliptic curves in families of quadratic twists \cite{Goldfeld}, states that asymptotically half of all curves have rank $0$ and half have rank $1$, so the average rank is exactly $1/2$.

The sign appearing in the functional equation for the $L$-function associated to an elliptic curve $E$ over $\Q$ is the {\em parity} of $E$, and it is conjectured that among all curves ordered by any reasonable invariant (like conductor, discriminant, or height), asymptotically half will have each parity.  The Parity Conjecture, a consequence of the BSD Conjecture, states that the sign of the functional equation is equal to the parity of the Mordell-Weil rank of $E$.  Therefore, conjecturally half of all curves have odd rank and thus have rank at least $1$, so the average rank of curves ordered by conductor is at least $1/2$. 

The Minimalist Conjecture also follows from this conjectural equidistribution of parity and the idea that elliptic curves generally have as small of a rank as allowed by parity; see also \cite{BMSW} for an extended discussion. A consequence of the Minimalist Conjecture is that zero percent of curves (asymptotically) should have rank at least $2$. We discuss heuristics for higher rank curves in \S \ref{sec:rank2}.

\subsection{Selmer groups}

Recent breakthroughs involving orbit parametrizations of genus one curves and the geometry of numbers have led to new unconditional bounds on average ranks of elliptic curves ordered by naive height. These rank bounds are consequences of results on Selmer groups of elliptic curves.

For each integer $n \geq 2$, the $n$-Selmer group $S_n(E)$ of an elliptic curve $E$ over $\Q$ fits into an exact sequence
\begin{equation}\label{Exact}
0\rightarrow E(\Q)/nE(\Q) \rightarrow S_n(E) \rightarrow \Sha(E)[n] \rightarrow 0,
\end{equation}
where $\Sha(E)[n]$ denotes the $n$-torsion subgroup of the Tate-Shafarevich group $\Sha(E)$ of $E$ over $\Q$.  If $p$ is a prime, then $S_p(E)$ is an elementary abelian $p$-group, whose dimension as an $\F_p$-vector space is called the {\em $p$-Selmer rank} of $E$.

\begin{thm}[(Bhargava--Shankar \cite{BhaSha2,BhaSha1,BhaSha4Sel,BhaSha5Sel})]\label{BSnSel}
When all elliptic curves $E/\Q$ are ordered by naive height, for $n \le 5$, the average size of $S_n(E)$ is $\sigma(n)$, the sum of the divisors of $n$.
\end{thm}

In \S \ref{sec:SelSha}, we discuss the $2$-Selmer group sizes for elliptic curves in our database, ordered by height.  We believe this is the first large-scale database of Selmer group information to compare to these theoretical results. The average size of $\Sel_2(E)$ for all curves of height at most $2.7\cdot 10^{10}$ is $2.6656$ and seems to be increasing towards the theoretical asymptotic average of $3$.

A consequence of the Selmer group result for $n = 5$ in Theorem \ref{BSnSel} is an upper bound on the average Mordell-Weil rank.
\begin{cor}[(Bhargava--Shankar \cite{BhaSha5Sel})]
When all elliptic curves over $\Q$ are ordered by height, their average rank is at most $0.885$.
\end{cor}

Our data, especially the samples at larger height (see \S \ref{sec:samples1}), suggest that in fact the average rank of elliptic curves is well below $0.885$. Extending Theorem \ref{BSnSel} to all $n$ would imply the Minimalist Conjecture, and such a generalization is supported by heuristics such as \cite[Conjecture 1.3]{BhaKanLenPooRai}, obtained by modeling the exact sequence \eqref{Exact} of $\Z_p$-modules via random maximal isotropic spaces.

\subsection{Other families}
Bhargava and Ho have adapted these Selmer group arguments to apply to families of curves with marked points.  We define the following family of elliptic curves:
\begin{equation*}
\Fs_1 :=  \{E : y^2 + a_3 y = x^3 + a_2 x^2 +a_4 x\ |\ a_2, a_3, a_4 \in \Z,\ \Delta_E \neq 0\},\\
\end{equation*}
where $\Delta_E$ denotes the discriminant of the elliptic curve $E$. For this family, there is a natural height function
  $H_1(E) := \max\{a_2^6, a_3^4, |a_4|^3 \}$
for $E \in \Fs_1$. In \cite{BH}, they show that if elliptic curves in $\Fs_1$ are ordered by height $H_1$, then the average size of the $2$-Selmer groups is bounded above by $6$, the average size of the $3$-Selmer groups is $12$, and the average rank is bounded by $13/6$.

We have created a database of all isomorphism classes of elliptic curves in $\Fs_1$ with height $H_1 \leq 10^8$ and computed the same invariants, such as rank and $2$-Selmer rank. Note that only $\num{693601}$ (approximately $19.3\%$) of these curves are in the main database. This database is discussed in more detail in \S \ref{sec:markedpointfamily}.

In \cite{BH}, several other families of elliptic curves with marked points are also studied with similar results on the average sizes of Selmer groups and bounds on average ranks. It would be interesting to create databases for these families to compare with the theoretical results.

\subsection{Other properties}

Our database also includes several other invariants of elliptic curves for which we can give results similar to those above, e.g., the number of curves with a given torsion subgroup and rank. Table \ref{table:variousproperties} gives the number of elliptic curves $E$ of naive height at most $2.7\cdot 10^{10}$ having certain interesting properties, as well as the proportion of different ranks. We also may compute averages related to each of these properties: for example, the average $2$-rank of $\Sha(E)$ for the curves in the main database is $0.23912$ and the average rank of all curves having complex multiplication in this database is $0.89848$.

We note that $19.99\%$ of the curves in our database have positive discriminant and that the rank distributions for curves of positive and negative discriminant appear to have different behavior: the average rank of curves in our database with $\Delta_E< 0$ is $0.88694$ and the average rank of curves with $\Delta_E > 0$ is $0.961245$.  We discuss some of these issues further in \S \ref{sec:avgrank}.

\subsection{Samples} \label{sec:samples1}

\begin{figure}[ht]
\caption{Average rank of elliptic curves ($\log_{10}$ scale)}\label{fig:averagerankwithsamples}
\fig{0.95}{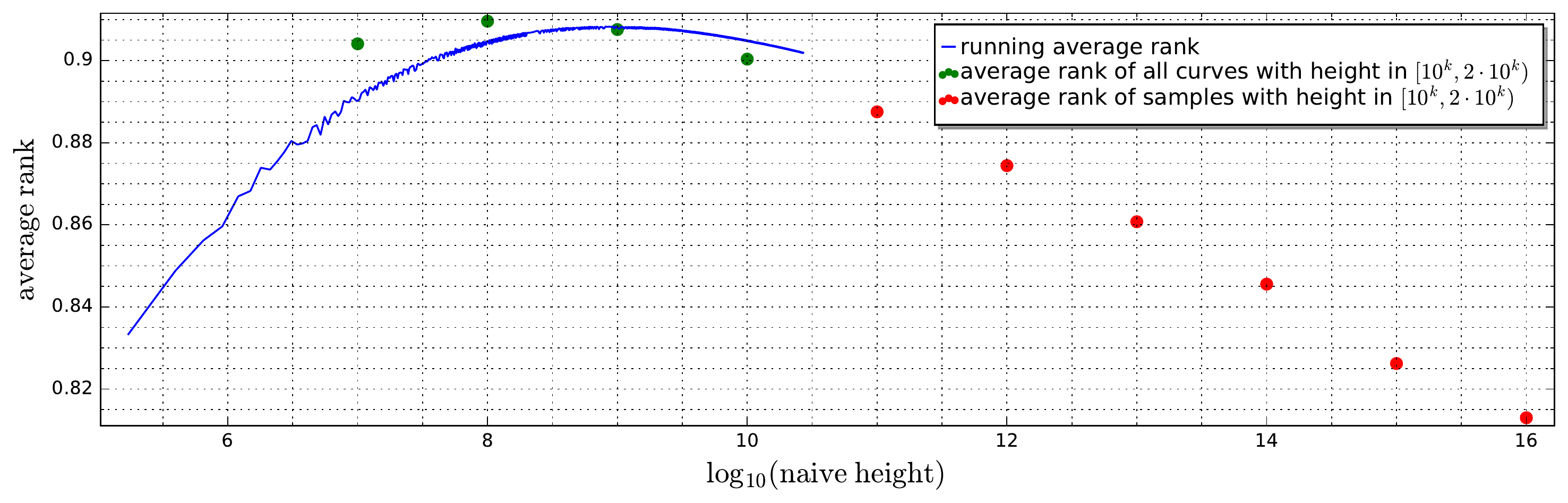}
\end{figure}

We have also created small databases of random samples of elliptic curves at larger heights. In particular, for each $k \in [11, 16]$, we chose $\num{100000}$ curves from a uniform distribution of all curves in the height range $[10^k,2 \cdot 10^k)$ and computed the same invariants.  These are not exhaustive datasets, but still provide some evidence for the behavior of various quantities as height increases. For example, we see the average rank decreases rapidly: see Table \ref{table:sampleranks} and Figure \ref{fig:averagerankwithsamples}, where the red points denote the average ranks for the samples, the green points denote the average rank for all curves in the height range $[10^k, 2 \cdot 10^k)$, and the blue curve represents the running average rank of all curves up to a given height.

We also record the distribution of the orders of the $2$-Selmer groups in Table \ref{table:sample2Selmerranks}; note the rapid convergence of the average $2$-Selmer size to the theoretical average of $3$. See Figures \ref{proprank2}, \ref{fig:averageselwithsamples}, and \ref{fig:averageshawithsamples} for the proportion of rank $2$ curves, the average $2$-Selmer size, and the average $2$-rank of $\Sha[2]$, respectively. Each figure includes the relevant values for these samples, denoted by red points; the green dots denote the corresponding values for all curves in the height range $[10^k, 2 \cdot 10^k)$ for $k = 7, 8, 9, 10$ and the blue curve represents the running average or proportion.

\section{Computing Ranks of Elliptic Curves}\label{Computing}

In this section, we describe our methods for populating our databases and computing information about Mordell-Weil groups and Selmer groups for the elliptic curves in the databases. There are two challenges involved when computing the rank: one must exhibit explicit rational points on $E$ while simultaneously showing that the rank of $E(\Q)$, denoted $\rk\,E(\Q)$, is no more than the rank of the subgroup generated by the known points.

There are two fundamentally different ways of obtaining upper bounds on the rank. The first relies on computing the $n$-Selmer group $S_n(E)$ for various integers $n$. This method gives the correct answer whenever it terminates, but its termination is conditional on the conjecture that $\Sha(E)$ is finite. The second relies on computing upper bounds on the order of vanishing of the $L$-function attached to $E$ and is conditional on the Birch and Swinnerton-Dyer conjecture.

\subsection{Populating databases}\label{Databases}

It is a straightforward task to write down all pairs of integers $(a_4,a_6)$ such that the naive height $\max\{4|a_4|^3, 27 a_6^2\}$ of the corresponding curve $y^2 = x^3 + a_4 x + a_6$ is in any chosen range.  For each such pair, we check that the discriminant of the corresponding curve does not vanish, i.e., the curve is nonsingular.  A curve of this form is isomorphic to one of smaller height if and only if there is a prime $p$ such that $p^4 \mid a_4$ and $p^6 \mid a_6$, and it is straightforward to remove non-minimal duplicates.  This process gives an exhaustive list of all isomorphism classes of curves in the desired height range.  For ease of computation and data analysis, we store the main database in approximately $30$ shards, most corresponding to a height range of size $10^9$.

We use a similar process to create an exhaustive list of isomorphism classes of curves in $\Fs_1$ of height at most $10^8$ using the modified height function $H_1$ for that family.

To create each of the larger height samples of $\num{100000}$ curves, for each integer $k \in [11,16]$, we repeatedly uniformly sample integers $a_4$ and $a_6$ from the appropriate ranges such that $4 |a_4|^3 < 2 \cdot 10^k$ and $27 a_6^2 < 2 \cdot 10^k$. If the curve $y^2 = x^3 + a_4 x + a_6$ is nonsingular, minimal, and has naive height at least $10^k$, then it is entered in the database.

\subsection{General procedure for computing rank}

The goal of this section is to explain how we compute the Mordell-Weil rank for each curve in our databases. We assume several conjectures during these computations---Birch and Swinnerton-Dyer (BSD), Generalized Riemann Hypothesis (GRH), and Parity---though not all of them are needed for every curve. Recall that the Parity Conjecture, which follows from BSD, states that for an elliptic curve $E$ over $\Q$ the root number is $(-1)^{\rk\,E(\Q)}$.

For each curve, we first compute standard arithmetic data, such as the conductor, root number, Tamagawa product, and torsion subgroup. To obtain rank, the first major step is to run Cremona's \texttt{mwrank} program with default parameters, which searches for points of low height and runs a $2$-descent. For each curve, \texttt{mwrank} yields the $2$-Selmer rank and upper and lower bounds for the Mordell-Weil rank.

If the bounds agree, we of course may determine the rank immediately, and if the difference between the upper and lower bound obtained from mwrank is $1$, then the root number combined with the Parity Conjecture gives the value of the rank. However, for many curves, the interval between the \texttt{mwrank} lower and upper bounds contains at least two integers of the ``correct'' parity, e.g., curves with even parity, lower bound $0$, and upper bound $2$.

In these cases, we attempt to improve the upper bound by applying the analytic technique described in \S \ref{sec:analytic}.  The upper bounds coming from this method are conditional on GRH.  In Corollary \ref{AnalyticRankUpperBound}, for any positive real parameter $\Delta$, we obtain an expression in terms of $\Delta$ that is an upper bound for the analytic rank of $E$, and which converges to the analytic rank from above.  Assuming BSD, the analytic rank is equal to $\rk\,E(\Q)$, so applying this bound with large enough $\Delta$ converges to the correct value of the rank.  Unfortunately, this method becomes computationally infeasible for large values of $\Delta$.  We compute this upper bound with successively larger values of $\Delta$, usually between $1$ and $3$, by the Sage function \texttt{analytic\_rank\_upper\_bound}, and stop the process and conclude that we have determined $\rk\,E(\Q)$ whenever the upper bound is within $1$ of the \texttt{mwrank} lower bound.  For a small number of curves in our larger height samples, we use this method with values of $\Delta$ up to $3.9$, which took several days for each curve at the highest values of $\Delta$.

This process allows us to conclude $\rk\,E(\Q)$ in the vast majority of cases.  For the remaining curves, we use methods in Magma to conclude the correct rank by either finding additional rational points to improve the lower bound, or computing the Cassels--Tate pairing between Selmer group elements to improve the upper bound. These techniques are described in \S \ref{sec:Magma}.

\enlargethispage{\baselineskip}

\subsection{Analytic upper bounds} \label{sec:analytic}

The analytic rank of an elliptic curve may be bounded from above by a certain explicit formula-derived sum over the nontrivial zeros of $\Les$, at the expense of having to assume GRH. We reproduce \cite[Lemma 2.1]{Bober}, which is a version of the explicit formula for elliptic curve $L$-functions akin to the Weil formulation of the Riemann-von Mangoldt explicit formula for $\zeta(s)$; a proof may be found in \cite[Theorem 5.12]{IwKo-2004}.
\begin{lemma}
\label{lem:exp_form_2}
Assume GRH.  Let $E$ be an elliptic curve over $\Q$, and let
\begin{equation}\label{def:bn}
b_n(E) = \begin{cases}
-\left(p^e + 1 - \#\widetilde{E}(\F_{p^e})\right)\cdot \log(p), & n=p^e\;\;\text{a prime power,} \\
0, & \text{otherwise.} \end{cases}
\end{equation}
where $\#\widetilde{E}(\F_{p^e})$ is the number of points on the (possibly singular) curve over the finite field of $p^e$ elements obtained by reducing $E$ modulo $p$.
Let $\gamma$ range over the imaginary parts of nontrivial zeros of $E$, and let $c_n = c_n(E) = \frac{b_n(E)}{n}$. Suppose that $f(z)$ is an entire function such that
\begin{itemize}
\item there exists a $\delta>0$ such that $f(x+iy) = O(x^{-(1+\delta)})$ for $|y|<1+\epsilon$ for some $\epsilon>0$, and
\item the Fourier transform of $f$, given by $\hat{f}(y) = \int_{-\infty}^{\infty} e^{-i x y}f(x)\; dx$, exists and is such that $\sum_{n=1}^{\infty} c_n \cdot \hat{f}\left(\log n\right)$ converges absolutely.
\end{itemize}
Then
\begin{equation*}
\sum_{\gamma} f(\gamma) = \frac{1}{\pi}\left[\log\left(\frac{\sqrt{N_E}}{2\pi}\right)\hat{f}(0) + \Re\int_{-\infty}^{\infty} \digamma(1+it)f(t) \; dt  + \frac{1}{2} \sum_{n=1}^{\infty} c_n \left( \hat{f}\left(\log n\right) + \hat{f}\left(-\log n\right)\right) \right],
\end{equation*}
where $\digamma(z)$ denotes the digamma function, the logarithmic derivative of $\Gamma(z)$.
\end{lemma}

We may use the above to provide computationally effective upper bounds on the analytic rank of an elliptic curve by choosing an appropriate test function $f$ whose Fourier transform has compact support. The method appears to have first been formulated by Mestre in \cite{Me-1986}, and used by Brumer in \cite{Bru-1992} to prove that, conditional on GRH, the average rank of elliptic curves is at most 2.3. The method was further refined to produce an upper bound of 2 by Heath-Brown in \cite{HeBr-2004} and then $25/14$ by Young in \cite{You-2006}.

Specifically, we use the parameterized Fej\'{e}r kernel as used by Mestre, Brumer, and Heath-Brown in the publications above and by Bober in \cite{Bober}:
\begin{equation} \label{eq:Fejer}
f_{\Delta}(x) = \sinc^2(\Delta x) = \left(\frac{\sin(\Delta \pi x)}{\Delta \pi x}\right)^2
\end{equation}
where $\Delta>0$ is the tightness parameter. Its Fourier transform is the triangular function
\begin{equation} \label{eq:FejerFourier}
\hat{f}_{\Delta}(y) = \begin{cases} \frac{1}{\Delta} \left(1 - \frac{|y|}{2\pi\Delta}\right) & |y|\le 2\pi\Delta, \\
0 & \text{otherwise.} \end{cases}
\end{equation}
Moreover, the integral $\Re\int_{-\infty}^{\infty} \digamma(1+it)f_{\Delta}(t) \; dt$ can be computed explicitly in terms of known constants and special functions:
\begin{equation} \label{eq:Redigamma}
\Re\int_{-\infty}^{\infty} \digamma(1+it) \cdot f_{\Delta}(t) \; dt = - \frac{\eta}{\pi\Delta} + \frac{1}{2\pi^2 \Delta^2}\left(\frac{\pi^2}{6} - \Li_2\left(e^{-2\pi \Delta}\right)\right).
\end{equation}
Here $\eta \approx 0.57722$ is the Euler-Mascheroni constant, and $\Li_2(x) = \sum_{n\ge 1} \frac{x^n}{n^2}$ is the logarithmic integral function.
Combining \eqref{eq:Fejer},\eqref{eq:FejerFourier}, and \eqref{eq:Redigamma} with Lemma \ref{lem:exp_form_2}, we obtain:
\begin{cor}\label{AnalyticRankUpperBound}
Assume GRH.  Let $\gamma$ range over the imaginary parts of the nontrivial zeros of $L_E(s)$, and let $\Delta > 0$. Then
\small
\begin{equation*}
\sum_{\gamma} \sinc^2(\Delta \gamma) = \frac{1}{\Delta \pi}\left[\left(-\eta + \log\left(\frac{\sqrt{N_E}}{2\pi}\right)\right)+ \frac{1}{2\pi \Delta}\left(\frac{\pi^2}{6} - \Li_2\left(e^{-2\pi \Delta}\right)\right)  + \sum_{n<e^{2\pi \Delta}} c_n \cdot \left(1-\frac{\log n}{2\pi \Delta}\right)\right]
\end{equation*}
\normalsize
and since $\sinc^2(0)=1$ and $\sinc^2(x)\to 0$ as $x \to \infty$, the sum converges to the analytic rank of $E$ from above as $\Delta \to \infty$.
\end{cor}

\begin{figure}[!ht]
    \fig{.9}{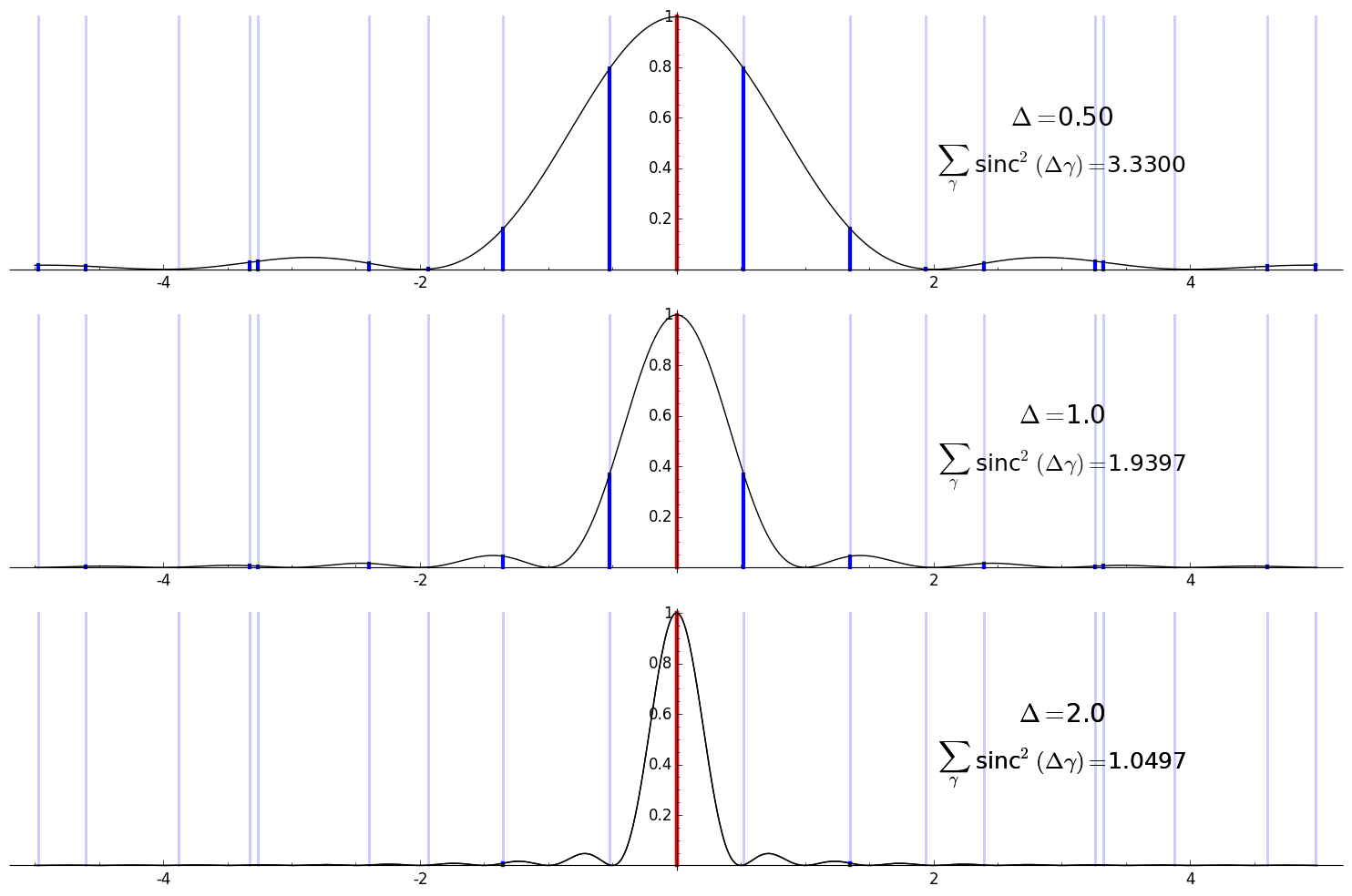}
    \caption{A graphic representation of the $\sinc^2$ sum for the elliptic curve $E: y^2=x^3-18x+51$, a rank $1$ curve with conductor $N_E=750384$, for three increasing values of the parameter $\Delta$. Vertical lines have been plotted at $x=\gamma$ whenever $L_E(1+i\gamma)=0$ (red for the single central zero and blue for noncentral zeros); the height of the darkened portion of each line is given by the black curve $\sinc^2(\Delta x)$. Summing up the lengths of the dark vertical lines thus gives the value of the $\sinc^2$ sum. As $\Delta$ increases, the contribution from the blue lines --- corresponding to noncentral zeros --- goes to zero, while the contribution from the central zero in red remains at 1. Thus the sum must approach 1 as $\Delta$ increases.}
    \label{fig:zero_sum_visualization}
\end{figure}

What is notable about the above formula is that evaluation of the right hand side is a finite computation and only requires knowledge of the elliptic curve's conductor and a finite number of $a_p$ values. See \cite[p. 67]{Spi-2015} for a more detailed derivation of this method. We may therefore (assuming BSD) obtain efficient upper bounds on the rank of an elliptic curve by choosing an appropriate value of $\Delta$ and performing this finite calculation.

\subsection{Magma techniques} \label{sec:Magma}

For a small number of curves in the main database and in the sample databases, we also use additional methods in Magma \cite{magma}.
While we assume GRH to speed up the construction of $2$-coverings and $4$-coverings of these curves, the rank bounds coming from these methods are unconditional, since we computed $\Sel_2(E)$ using \texttt{mwrank} and do not need to use a full list of $2$-coverings or $4$-coverings.

The Magma procedure begins by carrying out a 2-descent to compute the specific 2-coverings of $E$ corresponding to elements of $\Sel_2(E)$, and then searches for rational points of small height on the 2-coverings.

If the rank is not determined at this stage, then we compute the Cassels-Tate pairing on elements of $\Sel_2(E)/E(\Q)[2]$. Recall that the Cassels-Tate pairing $\Gamma$ is an alternating bilinear pairing on $\Sha(E)$ taking values in $\Q/\Z$; if $\Sha(E)$ is finite, then it is nondegenerate. When restricted to $\Sha(E)[2]$, this gives a nondegenerate alternating bilinear pairing on $\Sha(E)[2]/2\Sha(E)[4]$, or equivalently on $\Sel_2(E)/\mathrm{im}(\Sel_4(E))$, which takes values in $\Z/2\Z = \{0,1\}$. In particular, if $C$ and $D$ are 2-coverings of $E$ with $\Gamma(C, D) = 1$, then $C$ and $D$ correspond to elements of order 2 in $\Sha(E)$, which gives an improved upper bound on rank.

If necessary, we next use a $4$-descent to find explicit $4$-coverings of $E$ corresponding to some elements of $\Sel_4(E)$ and then search for rational points of small height on these $4$-coverings, refining lower bounds on the rank. For almost all curves, these methods, combined with the Parity Conjecture, are enough to determine the rank.

For $7$ curves in our sample database at height $10^{16}$, all these techniques, including computing analytic upper bounds with very large values of $\Delta$ and extensive point searches, do not determine rank: for each of these curves $E$, the $2$-rank of $\Sel_2(E)$ is $2$ and $E(\Q)$ has trivial torsion. We use Magma to compute the value $\Gamma(C,D) \in \Z / 2 \Z = \{0,1\}$ of the Cassels-Tate pairing for a $4$-covering $C \in S_4(E)$ and a $2$-covering $D \in S_2(E)$. We find an explicit pair $(C,D)$ with $\Gamma(C,D) = 1$, which implies that $\Sel_2(E) \cong \Sha(E)[2] \cong \Z/2\Z \times \Z/2\Z$ and thus $\rk\,E(\Q) = 0$. Two such curves for which this method is needed are $y^2 = x^3 + 169304 x + 25788938$ and $y^2 = x^3 + 77108 x -22146514$.


\section{Data Analysis}\label{Data Analysis}

\subsection{Distribution of rank} \label{sec:avgrank}

As mentioned in \S \ref{sec:intro}, we see from our main database that the average rank of elliptic curves of height up to $X$ increases with small $X$ and then decreases (on a large scale); see Figure \ref{average_rank_figure} for a plot of the average rank of all elliptic curves of height at most $2.7\cdot 10^{10}$ and Tables \ref{RankDistribution} and \ref{DatabaseAvg} for a more detailed distribution. Morever, in our larger height samples, the average rank decreases as height increases, with the average of the height $10^{16}$ sample approximately $0.813$ (see Table \ref{table:sampleranks} and Figure \ref{fig:averagerankwithsamples} for details).  For comparison, we note that the average rank of all curves of conductor up to $\num{360000}$ is $0.72759$ (using data from \cite{Cremona}) and the average rank of all curves of conductor at most $10^8$ in the Stein--Watkins database is $0.865$ \cite{BMSW}.

\subsubsection{Higher rank curves} \label{sec:rank2}

There are many proposed heuristics for predicting the asymptotics of rank $2$ curves. For example, when curves are ordered by conductor, Watkins \cite{Wat} predicts,  based on ideas from random matrix theory, that the number of isomorphism classes of elliptic curves with conductor at most $X$ and rank $2$ is $O(X^{19/24}(\ln X)^{3/8})$. For curves ordered by height, the recent heuristics of Park, Poonen, Voight, and Wood \cite{ParkPoonenVoightWood} predict that, for $2 \leq r \leq 20$, the number of curves with height up to $X$ and rank at least $r$ is asymptotic to $X^{(21-r)/24+o(1)}$ (but are not fine enough to predict the logarithmic term).

\begin{figure}[ht]
\caption{Proportion of rank 2 curves, including samples ($\log_{10}$ scale)}\label{proprank2}
\fig{.95}{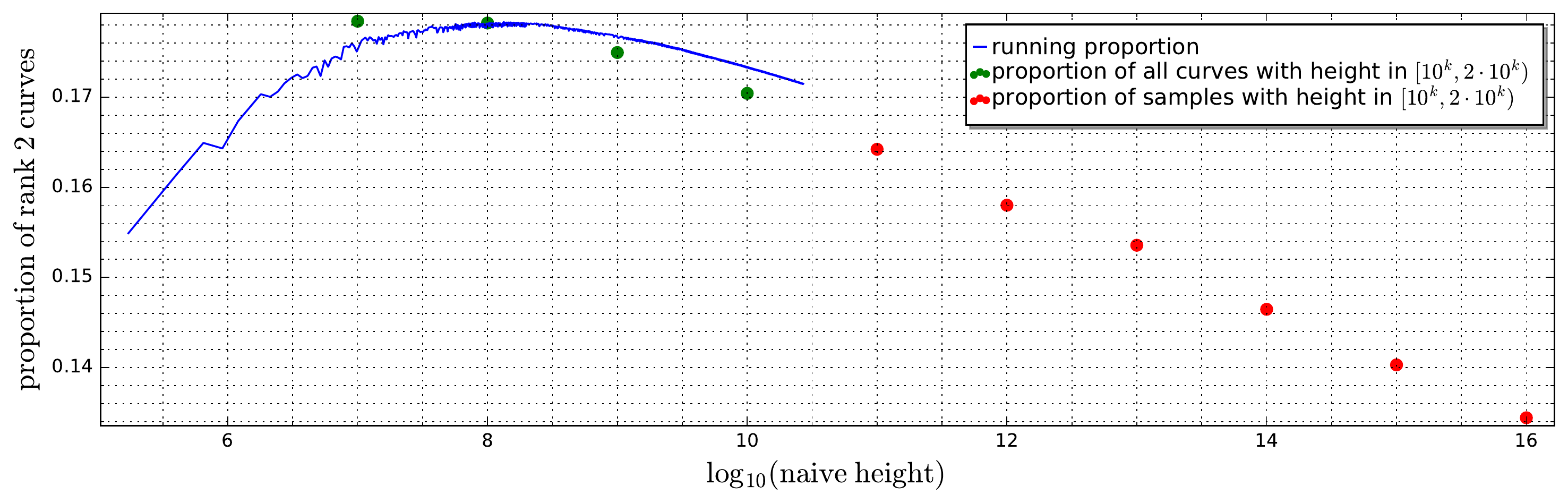}
\end{figure}

In our main database, there are $\num{40949307}$ rank $2$ curves, approximately $17.15\%$ of all of the curves (and $34.3\%$ of the even rank curves). The number of rank $2$ curves in the main database is approximately $0.0686 \cdot \left(2.7 \cdot 10^{10}\right)^{19/24} \left(\ln (2.7 \cdot 10^{10})\right)^{3/8}$, and we note that the constant seems to slightly increase as height increases. However, the proportion of curves of height up to $X$ having rank $2$ decreases as $X$ increases (for $X$ larger than approximately $10^8$): see Figure \ref{proprank2} and Table \ref{table:sampleranks}.

\subsubsection{Positive versus negative discriminant}

It is believed that asymptotically the distribution of ranks of elliptic curves $E$ with height at most $X$ and $\Delta_E>0$ should be the same as the distribution of ranks of elliptic curves $E$ with height at most $X$ and $\Delta_E<0$; however, for small values of $X$, these distributions initially appear different.  Brumer and McGuinness \cite{BrumerMcGuinness} note that in their database of $\num{310716}$ curves of prime conductor up to $10^8$, the average rank of those with $\Delta_E > 0$ is $1.04$, while the average rank for those with $\Delta_E <0$ is $0.94$. In \cite{BMSW}, the authors point out that the relationship between the sign of the discriminant and the average rank is a little subtle, computing far enough to find a crossing point in the graphs of average rank of curves of conductor at most $10^8$ in the Stein--Watkins database with given sign of $\Delta_E$.

\begin{figure}[ht]
\caption{Average rank of curves with positive discriminant versus negative discriminant ($\log_{10}$ scale)}\label{avgrankdisc}
\fig{.95}{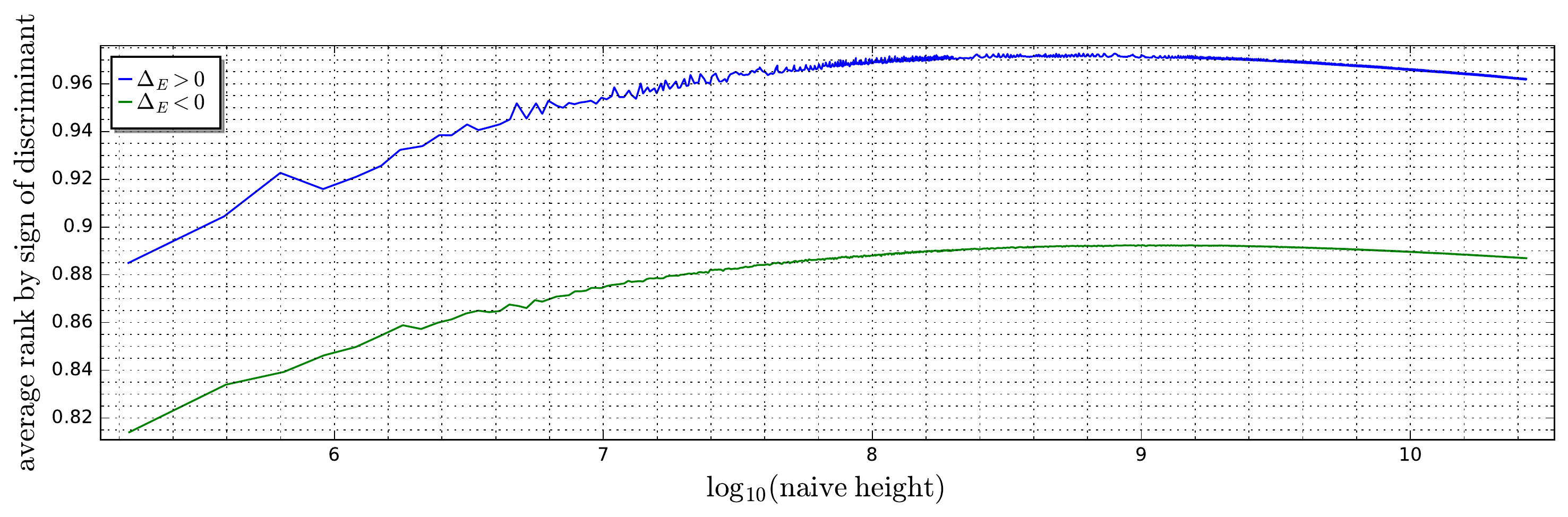}
\end{figure}

By the form of our height function, a curve $E$ has $\Delta_E > 0$ if and only if $4 |a_4|^3 > 27 a_6^2$ and $a_4<0$.  This accounts for exactly half of curves for which $4 |a_4|^3 > 27 a_6^2$, but less than half (in fact, $19.99\%$) of all curves in our main database.  Among all curves of naive height at most $2.7 \cdot 10^{10}$, we see that the average rank of those with $\Delta_E > 0$ is $0.961245$ while the average rank of those with $\Delta_E < 0$ is $0.88694$.  In fact, rank is weakly correlated with the sign of the discriminant, with a correlation coefficient $r=0.03856$; while this correlation value is small, it is still significant given the large size of the database. Note also that the fact that the sizes of these two sets of curves are not close to being equal does not explain this discrepancy in rank.  It would be interesting to have a theoretical explanation for these numerical observations.  Figure \ref{avgrankdisc} plots the average rank for all curves with each sign of $\Delta_E$ and height less than a given value.

\subsection{Selmer groups and Tate-Shafarevich groups} \label{sec:SelSha}

\begin{figure}[ht]
\caption{Average order of the $2$-Selmer groups, including samples ($\log_{10}$ scale)}\label{fig:averageselwithsamples}
\fig{.95}{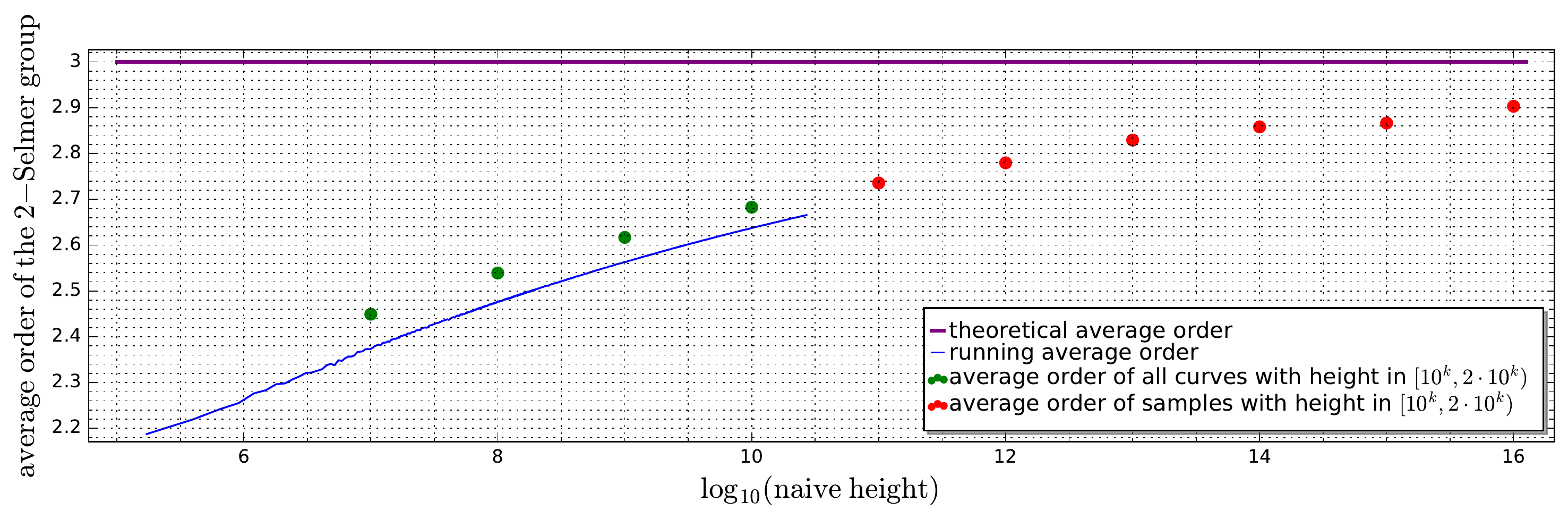}
\end{figure}

Theorem \ref{BSnSel} says that the average size of $\Sel_2(E)$ among all elliptic curves converges to $3$.  Since the rank of $\Sel_2(E)$ gives an upper bound on $\rk\,E(\Q)$ and the average rank of small height is larger than the conjectured asymptotic value of $1/2$, it may seem reasonable to guess that the average size of $\Sel_2(E)$ exceeds the theoretical average.  However, the average size of $\Sel_2(E)$ for all curves of height up to $X$ where $X \le 2.7 \cdot 10^{10}$ appears to be increasing towards the predicted value of $3$.  In our samples, the average size of $\Sel_2(E)$ increases in each larger height sample, with the average size in the $10^{16}$ sample already up to $2.90311$; see Table \ref{table:sample2Selmerranks} and Figure \ref{fig:averageselwithsamples}.

Conjecture 1.1 in \cite{PR} predicts that as we vary over all elliptic curves $E$ over $\Q$ ordered by height, we have
\begin{equation} \label{eq:PRprobSel}
\operatorname{Prob}\left(\dim_{\F_p} \Sel_p(E) = d\right) = \left( \prod_{j \ge 0} (1+p^{-j})^{-1}\right)\left(\prod_{j=1}^d \frac{p}{p^j-1}\right),
\end{equation}
which is compatible with the more general conjectures of \cite{BhaKanLenPooRai}.  For $p=2$, equation \eqref{eq:PRprobSel} predicts that the proportion of curves with $\dim_{\F_2} \Sel_2(E) = 0,1$, and $2$ should be approximately $0.2097,\ 0.4194$, and $0.2796$, respectively.  We see that the proportion of curves in our main database with $\dim_{\F_2} \Sel_2(E) = 0,1$, and $2$ are approximately $0.2381$, $0.4449$, and $0.2578$, respectively, quite close to the predicted values.

Since we record $\rk\,E(\Q)$, the size of $\Sel_2(E)$, and the torsion subgroup $E(\Q)_{\tors}$ for each elliptic curve $E$ in our database, we also easily deduce the size of the $2$-torsion part $\Sha(E)[2]$ of the Tate-Shafarevich group of $E$. Delaunay gives a conjecture for the asymptotic distribution of $\rk_{p^j}\Sha(E)$ in analogy with the Cohen-Lenstra heuristics for class groups of number fields \cite{Delaunay}.  More precisely, as we vary over all curves $E$ over $\Q$ up to isomorphism of rank $r$ ordered by conductor, he predicts
\[
\operatorname{Prob}\left(\dim_{\F_p} \Sha(E)[p] = 2n\right) = p^{-n(2r+2n-1)} \frac{\prod_{i=n+1}^\infty (1-p^{-(2r+2i-1)})}{\prod_{i=1}^n (1-p^{-2i})}.
\]
See \cite[Conjecture 5.1]{PR} for this version of the statement, where it is noted that it is reasonable to expect the same result to hold for curves ordered by height, or \cite[Conjecture 4]{DelJou} for a slightly different phrasing. For example, these heuristics predict that for rank $0$ curves, the proportion with $2$-rank of $\Sha(E)$ equal to $0$, $2$, or $4$ is equal to $0.4194$, $0.5592$, or $0.0213$, respectively, and for rank $1$ curves, the proportion with $2$-rank $0$ or $2$ is $0.8388$ or $0.1598$, respectively.  The moments of the conjectured distribution of $|\Sha(E)(p^j)|$ are then computed by Delaunay and Jouhet as \cite[Conjecture 3]{DelJou}: the expected value of $|\Sha(E)[2]|$ for curves of rank $r$ is $1+2^{-(2r-1)}$.

In our main database, the proportions of curves with rank $r = 0$ and $\dim_{\F_2} \Sha(E)[2] = 2n$ for $n = 0,1$, and $2$ are $0.7294$, $0.2695$, and $0.0011$, respectively, and with $r = 1$ and $n = 0,1$ are $0.9393$ and $0.0607$, respectively.  The average size of $|\Sha(E)[2]|$ for these rank $0$ curves is $1.825$ and for rank $1$ is $1.182$.  We see that the rank $0$ distribution of $\Sha(E)[2]$ is not particularly close to the theoretical predictions, but that the data fit more closely for curves of rank $1$.  In each case the average size of $\Sha(E)[2]$ is significantly smaller than expected, which helps to explain why the average size of $\Sel_2(E)$ appears to approach the asymptotic value of $3$ from below, even though the average rank seems to approach the asymptotic value of $1/2$ from above.

\begin{figure}[ht]
\caption{Average $2$-rank of $\Sha[2]$, including samples ($\log_{10}$ scale)}\label{fig:averageshawithsamples}
\fig{.95}{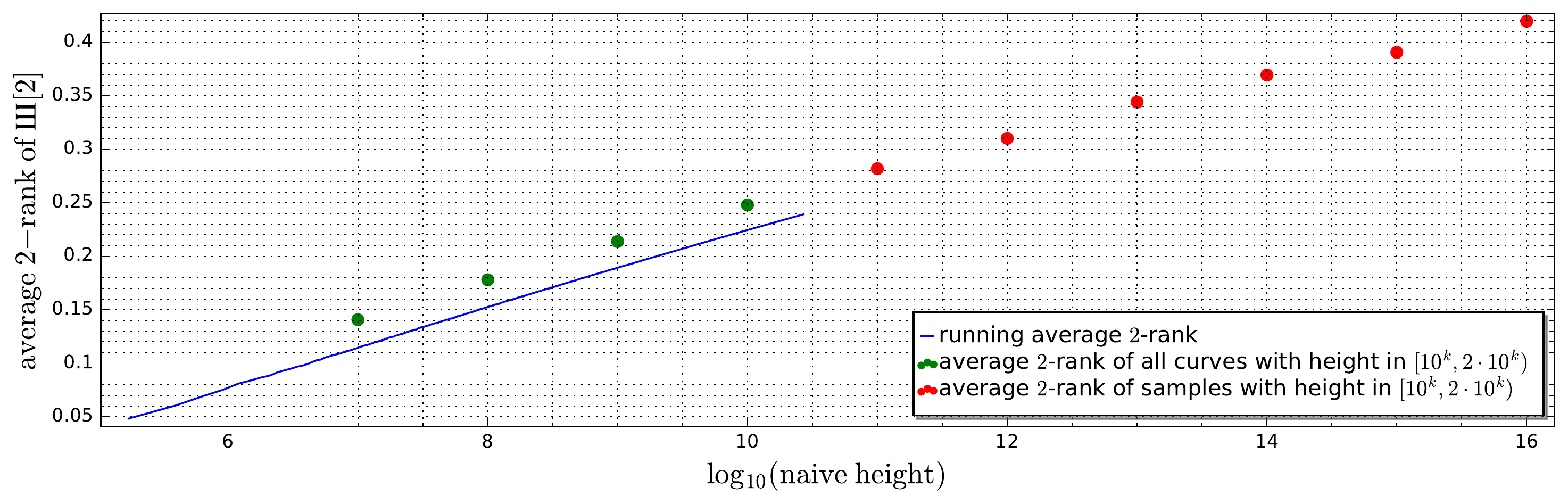}
\end{figure}

See Figure \ref{fig:averageshawithsamples} for a plot of the average $2$-rank of $\Sha[2]$ up to a given height, and Figure \ref{fig:averageshasizebyrrank} for a plot of the average size of $\Sha[2]$ up to a given height for rank $0$ and rank $1$ curves separately.

\subsection{Other invariants: root number and torsion subgroups}

One may ask whether the convergence of other arithmetic invariants appears to be faster than the convergence of average rank or average size of the $2$-Selmer group. For example, it is natural to conjecture that unless there is a good reason to believe that they must be biased, the root numbers of elliptic curves in families should be equidistributed. As noted earlier, in our main database the proportion of curves with root number $1$ is $0.49995$, already quite close to the conjectured value of $1/2$. Figure \ref{avgrootplot} shows how the average root number appears to quickly converge to the expected theoretical value of $0$.

Another example comes from studying the torsion subgroups of elliptic curves.  We know that as $X \rightarrow \infty$ the average size of the torsion subgroup of an elliptic curve of height up to $X$ approaches $1$; see Figure \ref{average_torsion} where this convergence appears to be quite fast.  We recall a more precise theorem:

\begin{thm}[(Harron-Snowden \cite{HarSno})]
Consider all elliptic curves in $\Fs_0$ ordered by uncalibrated height. Let $N_G(X)$ be equal to the number of isomorphism classes of elliptic curves of height up to $X$ with torsion subgroup $G$.  Then
\begin{equation*}
N_{\text{trivial}}(X) \sim \frac{4}{\zeta(10)} X^{5/6}, \qquad N_{\Z/2\Z}(X) \sim  c_2 X^{1/2}, \qquad \text{ and } \qquad
N_{\Z/3\Z}(X) \sim  c_3 X^{1/3},
\end{equation*}
where $c_2 \approx 3.1969$ and $c_3 \approx 1.5221$.
\end{thm}
We emphasize that this result uses the uncalibrated height function; our database includes all $\num{126427408}$ elliptic curves of uncalibrated height at most $10^9$. Table \ref{table:numberwithtorsion} includes the number and average of these elliptic curves with each possible torsion subgroup.

The number of curves with trivial torsion and uncalibrated height at most $10^9$ is approximately $0.9995\cdot \frac{4}{\zeta(10)} (10^9)^{5/6}$.  Similarly, the numbers with uncalibrated height at most $10^9$ and $E(\Q)_{\tors} \cong \Z/2\Z$ and with $E(\Q)_{\tors} \cong \Z/3\Z$ are approximately $1.2125 \cdot 3.1969 (10^9)^{1/2}$ and $0.4993 \cdot 1.5221 (10^9)^{1/3}$, respectively.  The number of curves with trivial torsion and bounded uncalibrated height appears to converge very quickly to the theoretical value, with the convergence being slower for torsion subgroups that occur less frequently in $\Fs_0$.

We note that if we restrict to curves with a given torsion subgroup, the analogue of the Minimalist Conjecture is expected to hold, implying that half of all curves have rank $0$ and half have rank $1$.  We find that the average rank of all curves of naive height at most $2.7\cdot 10^{10}$ and $E(\Q)_{\tors} = \Z/2\Z$ is $0.79895$ and the average rank of those with $E(\Q)_{\tors} = \Z/3\Z$ is $0.60882$. Figures \ref{avgrank2tor} and \ref{avgrank3tor} show plots of average rank of curves with naive height up to $X$ and torsion subgroups $\Z/2\Z$ and $\Z/3\Z$, respectively.

\subsection{Elliptic curves with complex multiplication} \label{sec:CM}

In \cite{BMSW}, the authors also consider average rank statistics for elliptic curves with complex multiplication (CM), those which have endomorphism ring (over $\C$) strictly larger than $\Z$. In the Stein--Watkins database of curves with conductor at most $10^8$, the proportion of the set of $\num{135226}$ curves with CM that have rank $0$ is $0.411$, significantly larger than the overall proportion $0.336$ of curves of rank $0$.  The average rank of the CM curves in that database is $0.687$.

In our main database of curves with naive height at most $2.7 \cdot 10^{10}$, there are $\num{65732}$ curves with CM, with only $32.819\%$ of them having rank $0$.  In fact, the rank distribution for the CM curves looks approximately like that of the entire main database, as expected, and the average rank of these CM curves is $0.89848$. Figure \ref{avgrankcm} gives a plot of average rank up to a given height for these CM curves.

\subsection{Family of elliptic curves with one marked point} \label{sec:markedpointfamily}

For the family $\Fs_1$ of elliptic curves with a marked point, the rank and $2$-Selmer rank distribution for the $\num{3594891}$ isomorphism classes of elliptic curves in $\Fs_1$ with height at most $10^{8}$ are as follows:
\[
\begin{tabular}{|c|c|c|}
Rank & No.~of Curves & \% of Curves \\
\hline
$0$ & $15783$ & $0.44\%$  \\
$1$ & $1239600$ & $34.48$ \\
$2$ & $1724209$ & $47.96$ \\
$3$ & $564784$ & $15.71$ \\
$4$ & $49642$ & $1.38$ \\
$5$ & $872$ & $0.024$ \\
$6$ & $1$ & $0.000028$ \\
\end{tabular}
\qquad
\begin{tabular}{|c|c|}
$2$-Selmer Rank & No.~of Curves \\
\hline
$0$ & $364$   \\
$1$ & $1145633$  \\
$2$ & $1727290$  \\
$3$ & $657323$  \\
$4$ & $63235$  \\
$5$ & $1045$  \\
$6$ & $1$  \\
\end{tabular}.
\]
The average rank of these curves is $1.83185$ and the average size of the $2$-Selmer group is $4.31296$.
Note that the $\num{15783}$ rank $0$ curves here all have nontrivial torsion, though the marked point is asymptotically a non-torsion point $100\%$ of the time by Hilbert irreducibility.

An analogue of the Minimalist Conjecture predicts that the average rank among all curves in $\Fs_1$ converges to $3/2$.  Just as in the family $\Fs_0$ in our main database, we see that the average rank of curves with ``small'' height is larger than the expected asymptotic value. Notably, despite having many fewer curves in this database than in our main database, the distribution here is closer to the asymptotic expectation, e.g., there are $15.71\%$ rank $3$ curves here, compared to $17.15\%$ rank $2$ curves in the main database.

\enlargethispage{\baselineskip}

\begin{acknowledgements}\label{ackref}
We used the open-source software SageMath \cite{Sage} and SageMathCloud \cite{SMC} extensively throughout this project, and we thank Harald Schilly for always keeping a close eye on our SMC projects. We are very grateful to the University of Michigan Advanced Research Computing Technology Services and the University of Oxford for providing further computing resources. We also thank Steve Donnelly for discussions about computing rank in Magma \cite{magma}. Finally, we thank the anonymous referees who provided many helpful suggestions for improving the paper.
\end{acknowledgements}

\clearpage

\setcounter{table}{0}
\makeatletter
\renewcommand*{\thetable}{A\arabic{table}}
\renewcommand*{\thefigure}{A\arabic{figure}}
\let\c@table\c@figure
\makeatother

\oneappendix 

\section{Additional tables and plots}

\begin{table}[ht]
\caption{Rank distribution for isomorphism classes of elliptic curves of naive height $\leq X$}
\label{RankDistribution}
\begin{tabular}{ccccccccc}
\hline
$X$ & Rank $0$ & Rank $1$ & Rank $2$ & Rank $3$ & Rank $4$ & Rank $5$ & Rank $6$\\
\hline
$10^8$ & $722275$ & $1073502$ & $400769$ & $51258$ & $1551$ & $7$ & $0$ \\
$10^9$ & $4930963$ & $7268430$ & $2706491$ & $384928$ & $16975$ & $137$ & $0$ \\
$10^{10}$ & $33944219$ & $49473528$ & $18099044$ & $2727260$ & $153537$ & $2119$ & $1$\\
$2 \cdot 10^{10}$ & $60667897$ & $88095239$ & $31992709$ & $4871438$ & $289954$ & $4654$ & $5$\\
$2.7 \cdot 10^{10}$ & $78039852$ & $113128980$ & $40949307$ & $6259159$ & $380519$ & $6481$ & $12$\\
\hline
\end{tabular}
\end{table}

\begin{table}[ht]
\caption{Average rank of isomorphism classes of elliptic curves of naive height $\leq X$}
\label{DatabaseAvg}
\begin{tabular}{|c|c|}
\hline
$X$ & Average rank of elliptic curves of naive height $\leq X$ \\
\hline
$10^8$ & $0.904724540$ \\
$10^9$ & $0.908338779$ \\
$10^{10}$ & $0.904965606$ \\
$2\cdot 10^{10}$ & $0.902949521$ \\
$2.7\cdot 10^{10}$ & $0.901975777$ \\
\hline
\end{tabular}
\end{table}

\begin{table}[ht]
\caption{Rank distribution in samples of \num{100000} elliptic curves of height between $10^k$ and $2 \cdot 10^k$}
\label{table:sampleranks}
\begin{tabular}{cccccccc}
\hline
$k$ &  Rank $0$ & Rank $1$ & Rank $2$ & Rank $3$ & Rank $4$ & Rank $5$ & Average rank\\
\hline
$11$ & $33318$ & $47547$ & $16422$ & $2495$ & $213$ & $5$ & $0.88753$ \\
$12$ & $34018$ & $47470$ & $15801$ & $2483$ & $219$ & $9$ &  $0.87442$ \\
$13$ & $34481$ & $47665$ & $15357$ & $2298$ & $192$ & $7$ & $0.86076$ \\
$14$ & $35000$ & $47991$ & $14647$ & $2180$ & $178$ & $4$ & $0.84557$ \\
$15$ & $35941$ & $47856$ & $14029$ & $1994$ & $174$ & $6$ & $0.82622$ \\
$16$ & $36407$ & $48105$ & $13442$ & $1885$ & $155$ & $6$ & $0.81294$ \\
\hline
\end{tabular}
\end{table}

\begin{table}[ht]
\caption{$2$-Selmer ranks in samples of $\num{100000}$ elliptic curves of height between $10^k$ and $2 \cdot 10^k$}
\label{table:sample2Selmerranks}
\begin{tabular}{cccccccc}
\hline
& \multicolumn{6}{c}{ --------------- $2$-rank of $2$-Selmer group ---------------} & Average size of\\
$k$ &  $0$ &  $1$ &  $2$ &  $3$ &  $4$ &  $5$ & $2$-Selmer group\\
\hline
$11$ & $23058$ & $44020$ & $26363$ & $6015$ & $532$ & $12$ & $2.73566$ \\
$12$ & $22829$ & $43541$ & $26608$ & $6392$ & $605$ & $25$ & $2.77959$ \\
$13$ & $22231$ & $43257$ & $27069$ & $6692$ & $729$ & $22$ & $2.82925$ \\
$14$ & $21973$ & $43177$ & $27073$ & $6968$ & $777$ & $32$ & $2.85819$ \\
$15$ & $22162$ & $42750$ & $27193$ & $7077$ & $786$ & $32$ & $2.86650$ \\
$16$ & $21613$ & $42631$ & $27553$ & $7329$ & $836$ & $38$ & $2.90311$ \\
\hline
\end{tabular}
\end{table}

\begin{table}[ht]
\caption{The number of elliptic curves $E$ in the main database with various properties and the proportion of curves with each rank out of those with the specified property.}
\label{table:variousproperties}
\begin{tabular}{|c|c|c|c|c|c||}
\hline
Property & No.~of Curves & Rank $0$ & Rank $1$ & Rank $2$ & Rank $\geq 3$ \\
& (\% of Database) \\
\hline
$E(\Q)_{\tors}$ trivial & $238528817$ & $0.327$ & $0.474$ & $0.172$ & $0.028$ \\
&($99.901\%$)\\
$E(\Q)_{\tors} \cong \Z/2\Z$ & $233153$ & $0.359$ & $0.492$ & $0.141$ & $0.008$ \\
& ($0.098\%$) \\
$E(\Q)_{\tors} \cong \Z/3\Z$ & $1020$ & $0.463$ & $0.466$ & $0.072$ & $0$\\
$E(\Q)_{\tors} \cong \Z/4\Z$ & $257$ & $0.521$ & $0.463$ & $0.016$ & $0$ \\
$E(\Q)_{\tors} \cong \Z/6\Z$ & $23$ & $0.870$ &	$0.130$ & $0$ & $0$ \\
$E(\Q)_{\tors} \cong \Z/2\Z \times \Z/2\Z$ & $1035$ & $0.453$ & $0.496$ & $0.051$ & $0$ \\
\hline
$\Delta_E > 0$ & $47738800$ & $0.305$ & $0.466$ & $0.192$ & $0.036$ \\
& ($19.994\%$)\\
$\Delta_E < 0$ & $191025510$ & $0.332$ & $0.476$ & $0.166$ & $0.026$ \\
& ($80.006\%$)\\
\hline
$\rk_2(\Sha(E)[2])=0$ & $210301413$ & $0.271$ & $0.505$ & $0.192$ & $0.032$ \\
& ($88.079\%$)\\
$\rk_2(\Sha(E)[2])=2$ & $28374370$ & $0.741$ & $0.242$ & $0.017$ & $0.00024$\\
& ($11.884\%$)\\
$\rk_2(\Sha(E)[2])=4$ & $88527$ & $0.978$ & $0.022$ & $0$ & $0$ \\
& ($0.037\%$) \\
\hline
$E$ has CM & $65732$ & $0.328$ & $0.474$ & $0.170$ & $0.028$ \\
\hline
$E$ has conductor $\leq 10^8$ & $4908673$ & $0.305$ & $0.474$ & $0.193$ & $0.027$ \\
\hline
\end{tabular}
\end{table}

\begin{table}[ht]
\caption{Elliptic curves $E$ of the form $y^2 = x^3 + a_4 x + a_6$ with minimal naive height for the specified rank and torsion subgroup.}
\label{Records}
\begin{tabular}{|c|c|c|c|c|c|}
\hline
$E(\Q)_{\tors}$ & Rank & $a_4$ & $a_6$ & $H(E)$ & Conductor($E$)\\
\hline
trivial & $0$ & $-1$ & $-1$ & $27$ & $368$ \\
trivial & $1$ & $-1$ & $1$ & $27$ & $92$ \\
&  & $1$ & $-1$ & $27$ & $248$ \\
&  & $1$ & $1$ & $27$ & $496$ \\
trivial & $2$ & $-4$ & $1$ & $256$ & $916$ \\
trivial & $3$ & $-13$ & $4$ & $8788$ & $66848$ \\
trivial & $4$ & $-19$ & $151$ & $615627$ & $4705528$ \\
trivial & $5$ & $-217$ & $1585$ & $67830075$ & $107827292$ \\
trivial & $6$ & $-1126$ & $6796$ & $5710513504$ & $35708014976$ \\
$\Z/2\Z$ & $0$ & $0$ & $1$ & $4$ & $64$ \\
$\Z/2\Z$ & $1$ & $-2$ & $0$ & $32$ & $256$ \\
$\Z/2\Z$ & $2$ & $7$ & $8$ & $1728$ & $4960$ \\
$\Z/2\Z$ & $3$ & $-82$ & $0$ & $2205472$ & $430336$\\
$\Z/2\Z$ & $4$ & $1030$ & $6396$ & $4370908000$ & $76983424$\\
$\Z/3\Z$ & $0$ & $0$ & $4$ & $432$ & $108$\\
$\Z/3\Z$ & $1$ & $0$ & $9$ & $2187$ & $972$\\
$\Z/3\Z$ & $2$ & $0$ & $225$ & $1366875$ & $24300$\\
$\Z/4\Z$ & $0$ & $-2$ & $1$ & $32$ & $40$\\
$\Z/4\Z$ & $1$ & $-2$ & $21$ & $11907$ & $760$\\
$\Z/4\Z$ & $2$ & $-191$ & $-510$ & $27871484$ & $7832$ \\
$\Z/5\Z$ & $0$ & $-43$ & $8208$ & $1819024128$ & $11$\\
$\Z/6\Z$ & $0$ & $0$ & $1$ & $27$ & $36$\\
$\Z/6\Z$ & $1$ & $-348$ & $2497$ & $168576768$ & $1260$ \\
$\Z/7\Z$ & $0$ & $-43$ & $166$ & $744012$ & $26$ \\
$\Z/9\Z$ & $0$ & $-219$ & $1654$ & $73864332$ & $54$ \\
$\Z/2\Z \times \Z/2\Z$ & $0$ & $-1$ & $0$ & $4$ & $32$ \\
$\Z/2\Z \times \Z/2\Z$ & $1$ & $-21$ & $-20$ & $37044$ & $288$ \\
$\Z/2\Z \times \Z/2\Z$ & $2$ & $-73$ & $72$ & $1556068$ & $19040$ \\
$\Z/2\Z \times \Z/4\Z$ & $0$ & $-351$ & $1890$ & $172974204$ & $24$  \\
\hline
\end{tabular}
\end{table}

\begin{table}[ht]
\caption{The number of elliptic curves $E$ in $\Fs_0$ with uncalibrated height at most $10^9$ and specified torsion subgroup}
\label{table:numberwithtorsion}
\begin{tabular}{|c|c|c|}
\hline
$E(\Q)_{\tors}$ & Number of Curves & Average Rank of these Curves \\
\hline
trivial & $126303317$ &  $0.894838$ \\
$\Z/2\Z$ &$122574$ &  $0.7832$ \\
$\Z/3\Z$ & $760$ & $0.59079$ \\
$\Z/4\Z$ & $188$ & $0.48936$ \\
$\Z/5\Z$ & $1$ & $0$ \\
$\Z/6\Z$ & $16$ & $0.125$\\
$\Z/7\Z$ & $1$ & $0$\\
$\Z/8\Z$ & $0$ & \\
$\Z/9\Z$ & $1$ & $0$ \\
$\Z/10\Z$ & $0$ & \\
$\Z/12\Z$ & $0$ & \\
$\Z/2\Z \times \Z/2\Z$ & $549$ &  $0.56466$ \\
$\Z/2\Z \times \Z/4\Z$ & $1$ & $0$  \\
$\Z/2\Z \times \Z/6\Z$ & $0$ &  \\
$\Z/2\Z \times \Z/8\Z$ & $0$ & \\
\hline
\end{tabular}
\end{table}

\clearpage
\begin{figure}[tp]
\caption{Average size of $\Sha[2]$ for curves of rank $0$ and rank $1$, including samples ($\log_{10}$ scale)}\label{fig:averageshasizebyrrank}
\fig{.95}{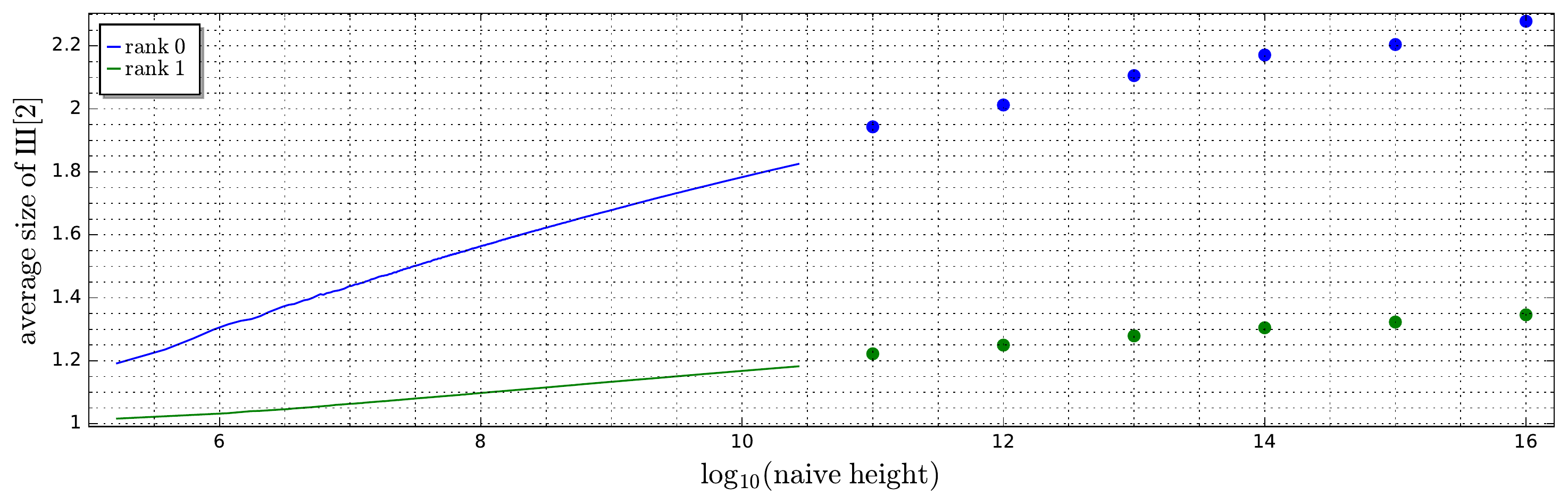}
\end{figure}

\begin{figure}[pptp]
\caption{Average root number ($\log_{10}$ scale)}\label{avgrootplot}
\fig{.95}{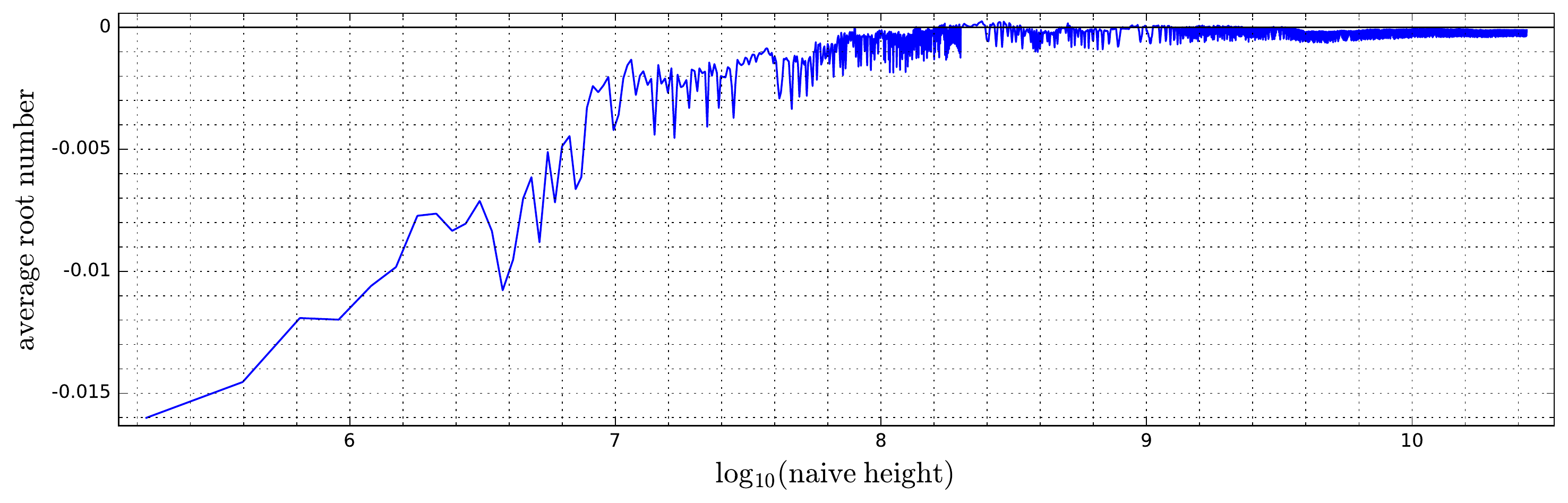}
\end{figure}

\begin{figure}[tppp]
\caption{Average order of the torsion subgroup ($\log_{10}$ scale)}\label{average_torsion}
\fig{.95}{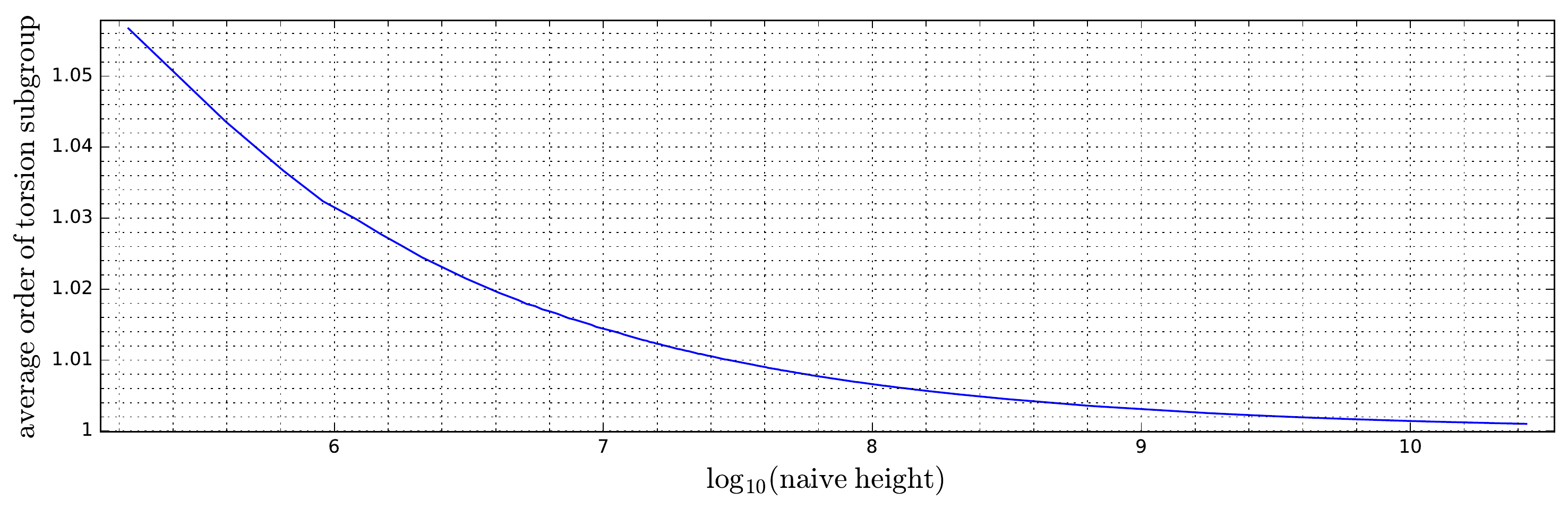}
\end{figure}

\begin{figure}[pptp]
\caption{Average rank of curves with torsion subgroup $\Z/2\Z$ ($\log_{10}$ scale)}\label{avgrank2tor}
\fig{.95}{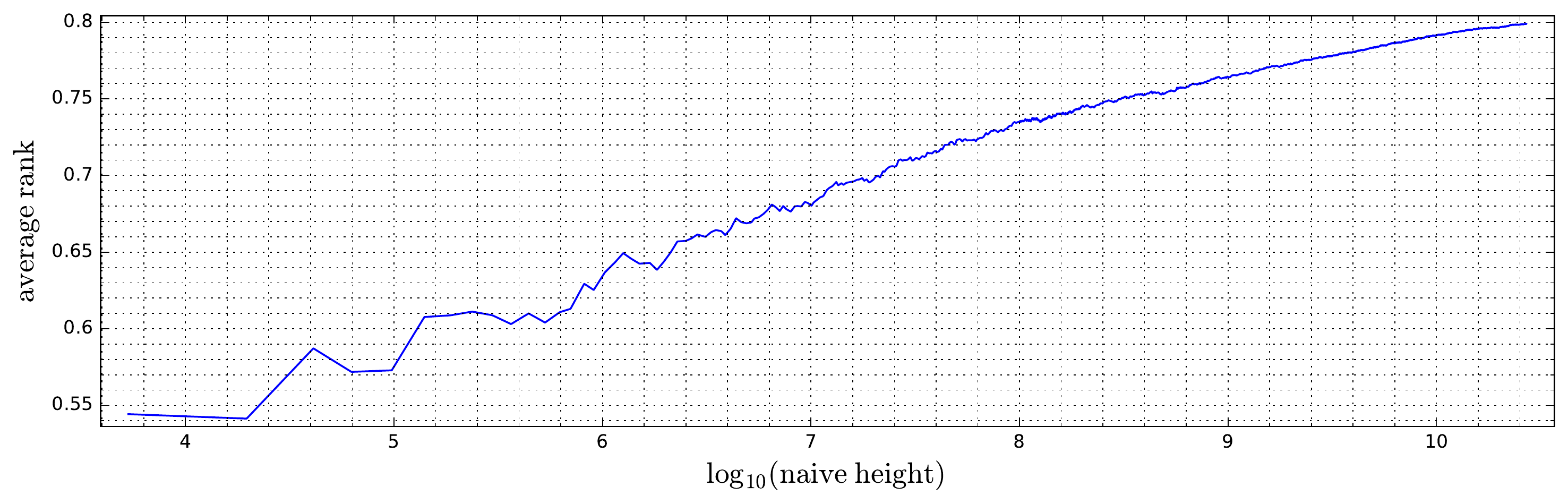}
\end{figure}

\begin{figure}[ht]
\caption{Average rank of curves with torsion subgroup $\Z/3\Z$ ($\log_{10}$ scale)}\label{avgrank3tor}
\fig{.95}{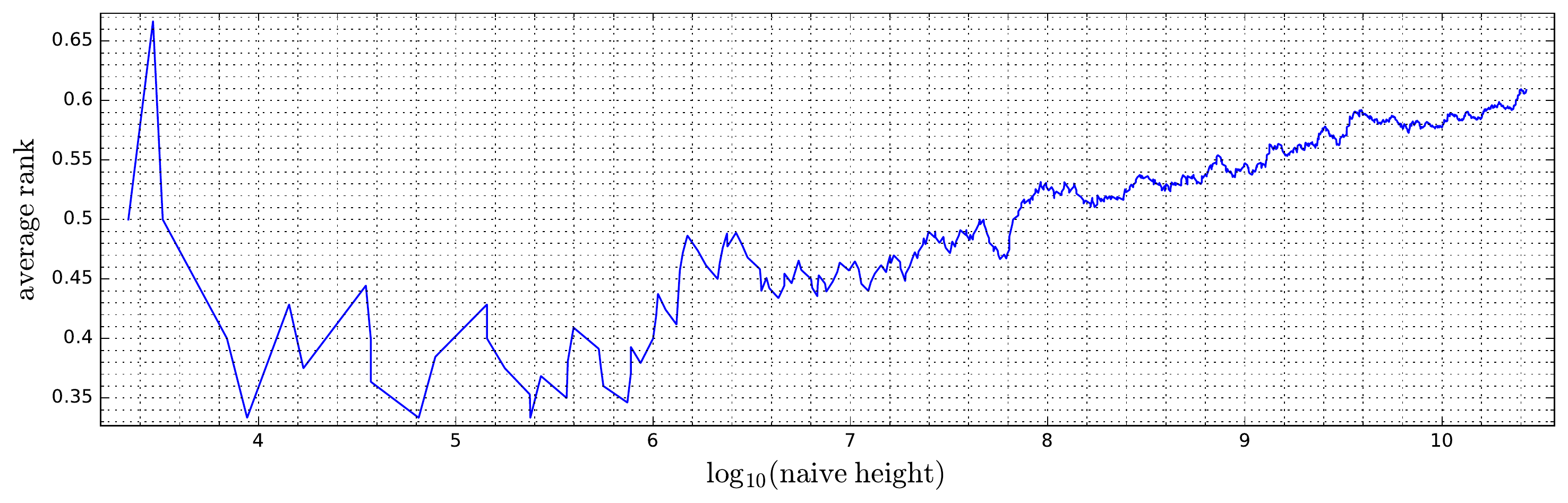}
\end{figure}

\begin{figure}[ht]
\caption{Average rank of CM curves ($\log_{10}$ scale)}\label{avgrankcm}
\fig{.95}{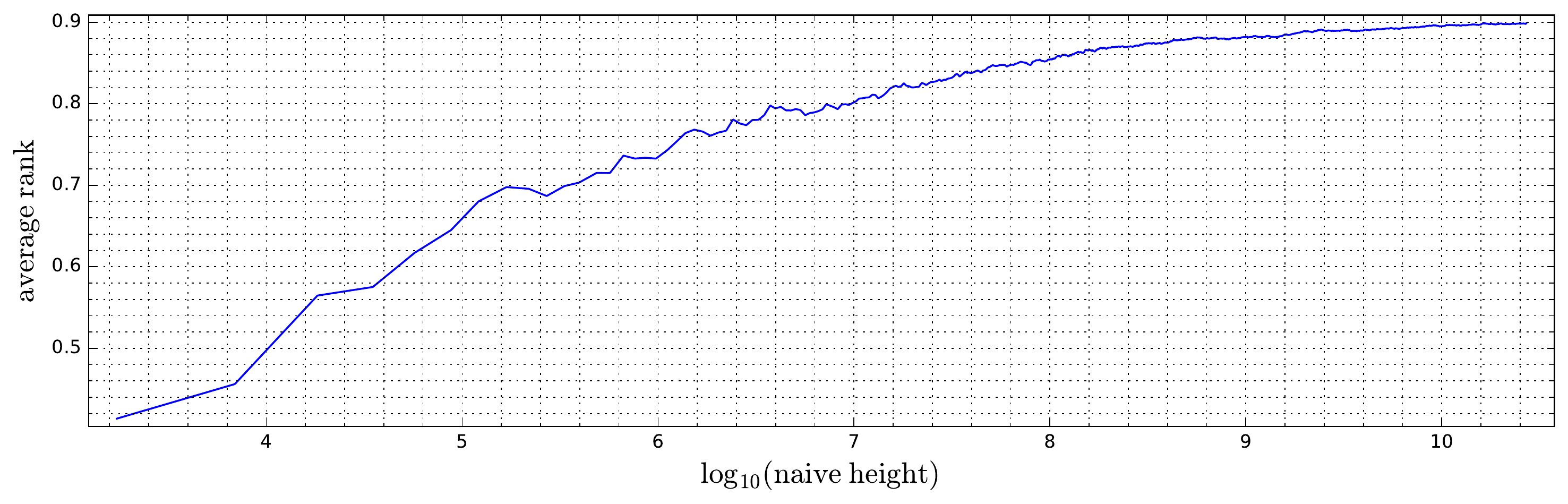}
\end{figure}

\begin{figure}[ht]
\caption{Average Tamagawa number ($\log_{10}$ scale)}\label{avgtam}
\fig{.95}{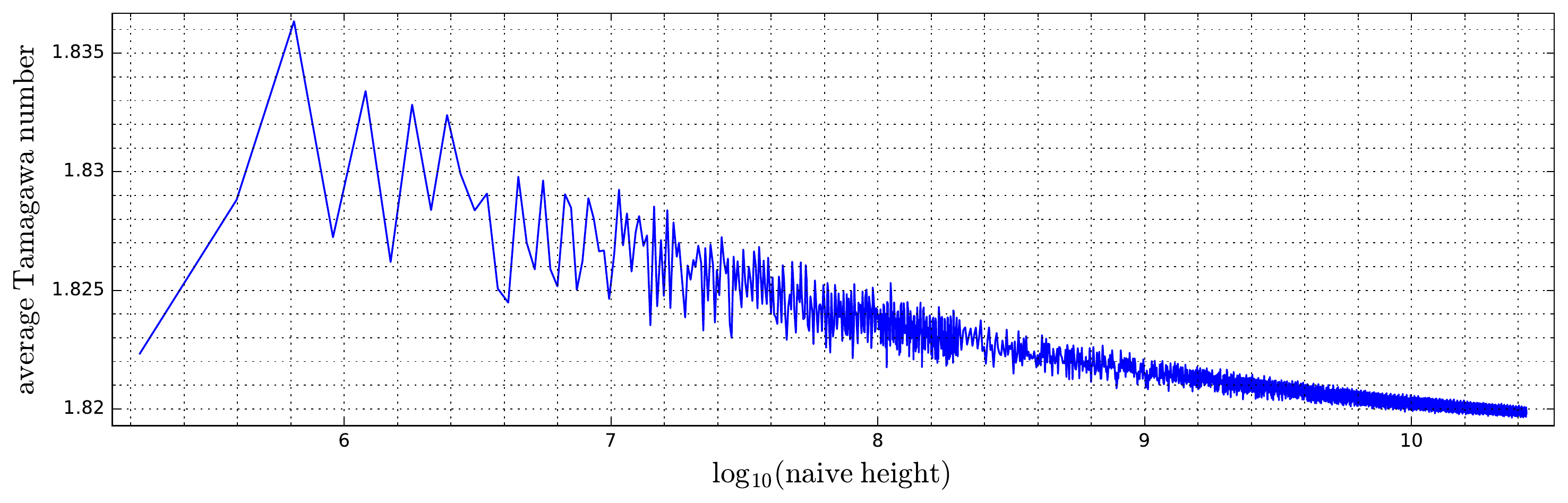}
\end{figure}

\bibliographystyle{amsplain}

\clearpage
\bibliography{bibliography}

\affiliationone{Jennifer S. Balakrishnan\\
Mathematical Institute\\
University of Oxford\\
Woodstock Road\\
Oxford OX2 6GG\\
United Kingdom
 \email{balakrishnan@maths.ox.ac.uk}}
\affiliationtwo{Wei Ho\\
Department of Mathematics\\
University of Michigan\\
Ann Arbor, MI 48109
USA
\email{weiho@umich.edu}}
\affiliationthree{Nathan Kaplan\\
Department of Mathematics\\
University of California, Irvine\\
Irvine, CA 92697
USA
\email{nckaplan@math.uci.edu}}
\affiliationfour{Simon Spicer\\
Facebook Inc.\\
1 Hacker Way \\
Menlo Park, CA 94025
USA
\email{mlungu@fb.com}}
\affiliationthree{William Stein\\
Department of Mathematics\\
University of Washington\\
Seattle, WA 98195-4350
USA
\email{wstein@uw.edu}}
\affiliationfour{James E. Weigandt\\
Institute for Computational and Experimental Research in Mathematics\\
Brown University\\
Providence, RI 02912
USA
\email{james\_weigandt@brown.edu}}

\end{document}